\documentclass[12pt]{article}
\usepackage{amssymb,amsmath}
\usepackage[all]{xy}
\newtheorem{thm}{Theorem}
\newtheorem{conj}{Conjecture}
\DeclareMathOperator{\des}{\mathbf{des}}
\DeclareMathOperator{\ides}{\mathbf{ides}}
\DeclareMathOperator{\run}{\mathbf{runs}}
\DeclareMathOperator{\asc}{\mathbf{asc}}
\DeclareMathOperator{\exc}{\mathbf{exc}}
\DeclareMathOperator{\inv}{\mathbf{inv}}

\begin{document}
\begin{center}
\Large
% TITLE GOES HERE
Two-sided Eulerian numbers via balls in boxes
\end{center}

\begin{flushright}
T. Kyle Petersen  \\
DePaul University \\
Chicago, IL 60614\\
\verb+tpeter21@depaul.edu+
\end{flushright}

\begin{quote}
I will shamelessly tell you what my bottom line is. It is placing balls into boxes\ldots \\

Gian-Carlo Rota, \emph{Indiscrete Thoughts}
\end{quote}

Gian-Carlo Rota was a professor at MIT from 1959 until his death in 1999. He is arguably the father of the field today known as algebraic combinatorics. Rota had 50 students (notably Richard Stanley) and, as of this writing, he has 329 mathematical descendants, including many of the top names in the field today. (And me.) If putting balls in boxes was good enough for him, well, it should be good enough for any of us. 

In this article we will first give a gentle introduction to permutation statistics and Eulerian numbers, followed by an explanation of how to use the idea of balls in boxes in this setting. Then we will take those ideas and study a perfectly natural (though less well-known) refinement: the ``two-sided" Eulerian numbers. We will finish with a discussion of symmetries of the Eulerian numbers and a conjecture of Ira Gessel.

\subsection*{Permutation statistics}

The symmetric group $S_n$ is the set of all permutations of length $n$, i.e., all bijections $w: \{1,2,\ldots,n\} \to \{1,2,\ldots,n\}$. We generally write a permutation in one-line notation: $w = w(1)w(2)\cdots w(n)$, so a typical element of $S_7$ is $w = 5624713$. 

For any permutation $w \in S_n$, we define a \emph{descent} to be a position $r$ such that $w(r) > w(r+1)$, and we denote by $\des(w)$ the number of descents of $w$. For example, if $w = 5624713$, then  there are descents in position 2 (since $6>2$) and in position 5 (since $7>1$). Hence, $\des(w) = 2$. The permutation $12\cdots n$ is the only permutation with no descents, while its reversal, $n \cdots 21$, has the maximal number of descents, with $n-1$.

The function \[\des: \quad \bigcup_{n \geq 1} S_n \to \{0,1,2,3,\ldots\} \] is an example of a \emph{permutation statistic}. Any function from the set of permutations to a subset of the integers can fairly be called a permutation statistic. The reader is referred to \cite{Bona} or \cite{K} for more about permutation statistics generally. 

Another common permutation statistic is \emph{inversion number}, $\inv(w)$, which is the number of pairs $r < s$ such that $w(r) > w(s)$. Notice the useful fact that $w$ and its inverse permutation must have the same number of inversions: if $r < s$ and $w(r) > w(s)$, then $w(s) < w(r)$ and $w^{-1}(w(s)) = s > r = w^{-1}(w(r))$. For example, with $w = 5624713$, we get $w^{-1}=6374125$ and we have $\inv(w)=\inv(w^{-1}) = 13$. (Notice that $\des(w)= 2$ while $\des(w^{-1})=3$, so descent numbers do not enjoy this property.) Again, the permutation $12\cdots n$ has no inversions, while $n\cdots 21$ has the most, with $\binom{n}{2}$ of them. Here, every pair of indices is an inversion. 

We mention $\inv$ because we can use it to give an alternate definition of $\des$ as follows. First, recall the \emph{$r$th simple transposition}, denoted $\sigma_r$, is the permutation that swaps $r$ and $r+1$ for some $r$ and fixes all other elements of $\{1,2,\ldots,n\}$. There are $n-1$ such elements in $S_n$, one for each $r=1,2,\ldots, n-1$. The simple transpositions are very special because they form a \emph{minimal generating set} for $S_n$, i.e., every permutation $w$ in $S_n$ can be obtained by applying a sequence of simple transpositions, and no proper subset of them will do the same.

Now consider the permutation obtained by composing a simple transposition with a permutation $w$: \[w\circ\sigma_r = w(1)\cdots w(r-1)w(r+1)w(r)w(r+2) \cdots w(n). \] We see this action swaps the numbers in positions $r$ and $r+1$ of $w$. Thus $r$ is a descent position of $w$ if and only if $r$ is \emph{not} a descent position of $w\circ\sigma_r$. In particular, $r$ is a descent of $w$ if and only if $\inv(w) > \inv(w\circ\sigma_r)$, i.e., \[ \des(w) = |\{ r : \inv(w\circ\sigma_r) < \inv(w) \}|.\] 

This identity is certainly not the most useful way to compute the descent number of a permutation, but it suggests that maybe we should consider the action of composition of $w$ with $\sigma_r$ on the left side as well. Here, we find $\sigma_r\circ w = (w^{-1}\circ\sigma_r)^{-1}$ since $\sigma_r^{-1}=\sigma_r$, and so (remembering inverses have the same inversion number) we have: \begin{align*}
 \des(w^{-1}) &= |\{ r : \inv(w^{-1}\circ\sigma_r) < \inv(w^{-1}) \}|\\
 &=|\{ r : \inv( \sigma_r \circ w) < \inv(w)\}|.
\end{align*}

Thus we may fairly call $\des(w)$ the number of \emph{right descents} of $w$, while $\des(w^{-1})$ is the number of \emph{left descents} of $w$. (The left descents are also called \emph{inverse descents} in the literature, and this number is also denoted by $\ides(w)$. See, for example, \cite[Section 13]{FH}.) In this article we will study both the right descent number on its own, and how left and right descents are related. This leads, respectively, to the Eulerian numbers, and to the ``two-sided" Eulerian numbers of the title.

\subsection*{The Eulerian numbers}

With a little patience, or a bit of computer programming, one can compute the distribution of $\des$ over $S_n$ by calculating $\des(w)$ for each $w$ in $S_n$. The reader may like to check by hand that in $S_3$ there is one permutation with no descents, four permutations with one descent, and one permutation with two descents. With a little more effort, you can find the distribution for $S_4$ as well. As this task requires $n!$ things to check for a given $n$, you may not want to try anything bigger than $n=4$ by hand.

We write $A_{n,i}$ for the number of permutations $w$ in $S_n$ for which $\des(w) = i-1$. (This shift in indices is done for historical reasons; primarily it makes certain formulas come out nicer.) The numbers $A_{n,i}$ are called the \emph{Eulerian numbers}, after Leonhard Euler. In Table \ref{tab:Eul}, we see the $A_{n,i}$ for $n\leq 8$.

\begin{table}[h]
\begin{center}
\begin{tabular}{r | c c c c c c c c c c c c}
$n \backslash i$ & 1 & 2 & 3 & 4 & 5 & 6 & 7 & 8\\
\hline
1 & 1 &&&&& \\
2 & 1 & 1 &&&&&\\
3 & 1 & 4 & 1 &&\\
4 & 1 &  11 & 11 & 1 \\
5 & 1 & 26 & 66 & 26 & 1\\
6 & 1 & 57 & 302 & 302 & 57 & 1\\
7 & 1 & 120 & 1191 & 2416 & 1191 & 120 & 1\\
8 & 1 & 247 & 4293 & 15619 & 15619 & 4293 & 247 & 1
\end{tabular}
\medskip
\caption{The Eulerian numbers $A_{n,i}$}\label{tab:Eul}
\end{center}
\end{table}

There is a great deal of mathematical literature devoted to the distribution of $\des$, and we remark that this distribution manifests itself in many other ways. For instance, let $\run(w)$ denote the number of \emph{increasing runs} of $w$ (maximal increasing subwords), let $\asc(w)$ denote the number of \emph{ascents} of $w$ (positions $r$ with $w(r)< w(r+1)$), and let $\exc(w)$ denote the number of \emph{excedances} of $w$ (positions $r$ with $w(r) > r$). Then, remarkably, we have:
\begin{align*}
A_{n,i} &=|\{ w \in S_n : \des(w) = i-1\}|,\\
 &=|\{ w \in S_n : \run(w) = i\}|,\\
 &=|\{ w \in S_n : \asc(w) = i-1\}|,\\
 &=|\{ w \in S_n : \exc(w) = i-1\}|.
\end{align*}
The second and third equalities are perhaps straightforward (there is necessarily one descent between each increasing run, and swapping each $w(r)$ for $n+1-w(r)$ trades descents for ascents), but that $\exc$ has the same distribution as $\des$ is not so obvious.

Euler himself was not interested in permutation statistics, but was rather investigating solutions of certain functional equations in which the numbers $A_{n,i}$ emerged. See the original \cite[Caput VII, pp. 389--390]{Eul} (available digitally online) or Carlitz's survey article, \cite{Car}. Carlitz points to several interesting occurrences of the Eulerian numbers in analysis and number theory in the first half of the 20th century, while in the early 1980s, Bj\"orner described the Eulerian numbers in a topological setting; see \cite{BjBr}. Suffice it to say, there are a variety of good reasons to study the Eulerian numbers, but for this article we will stick to the combinatorial point of view.

We will now discuss one of the standard tools of enumerative combinatorics, with which we will encode the Eulerian numbers. Given a sequence of numbers, $a_0, a_1, a_2, \ldots$, we define the \emph{generating function} for the sequence to be the power series (a polynomial, if the sequence is finite) in which the coefficient of $t^i$ is $a_i$: \[ \sum_{i\geq 0} a_i t^i = a_0 + a_1 t + a_2 t^2 + \cdots.\] So, for example, the geometric series \[ \frac{1}{1-t} = 1 + t + t^2 + \cdots, \] is the generating function for the sequence $1,1,1,\ldots$. Another example that might be familiar comes from the binomial theorem: \[ (1+t)^n = \binom{n}{0} + \binom{n}{1}t + \binom{n}{2}t^2 + \cdots + \binom{n}{n}t^n,\] so that for fixed $n$, $(1+t)^n$ is the generating function for the binomial coefficients $\binom{n}{i}$. Whether a generating function has finitely many terms or not, the reader should try not to worry about problems of convergence of the series. (For those who want to worry about it, rest assured that all series mentioned in this paper converge for $|t|<1$.) The reader should think of a generating function as a ``clothesline" on which the sequence hangs. 

The generating function for Eulerian numbers, with $n$ fixed, will be denoted \[A_n(t) = \sum_{i = 1}^n A_{n,i} t^i = A_{n,1}t + A_{n,2}t^2 + \cdots+ A_{n,n}t^n.\] We call this the \emph{$n$th Eulerian polynomial}. For example, the fourth Eulerian polynomial is: \[A_4(t) = t + 11t^2 + 11t^3 +t^4.\] Notice that since $A_{n,i}$ is the number of permutations with $i$ descents, we can get the same generating function by summing not on $i$, but on the set of all permutations in $S_n$, i.e., \[ A_n(t) = \sum_{w \in S_n} t^{\des(w)+1},\] since the number of summands contributing $t^i$ to the sum is precisely the number of permutations with $i-1$ descents.

We will now discuss the ``balls in boxes" approach to the Eulerian numbers.

\subsection*{Balls in boxes}

We begin with the generating function for the sequence $a_0, a_1, a_2, \ldots$, where $a_k$ is the number of ways of putting $n$ labeled balls into $k$ distinct boxes. 

We claim the generating function for all the ways of putting $n$ distinct balls into boxes is \[ a_0+a_1t + a_2t^2 + \cdots = \sum_{k\geq 0} k^n t^k.\]  Indeed, if there are $k$ different boxes, there are a total $k^n$ ways to put $n$ labeled balls into the boxes; we have $k$ choices for each ball. A typical choice with $k=9$, $n=6$ might look like this: 
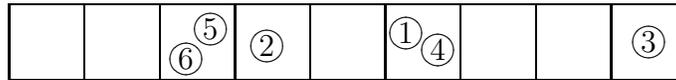
\[ 
\begin{xy}0;<.5cm,0cm>:
(-1,-1); (-1,1) **@{-}, (1,-1) **@{-}, (1,1); (1,-1) **@{-}, (-1,1) **@{-}
\end{xy}
\begin{xy}0;<.5cm,0cm>:
(-1,-1); (-1,1) **@{-}, (1,-1) **@{-}, (1,1); (1,-1) **@{-}, (-1,1) **@{-}
\end{xy}
\begin{xy}0;<.5cm,0cm>:
(-1,-1); (-1,1) **@{-}, (1,-1) **@{-}, (1,1); (1,-1) **@{-}, (-1,1) **@{-}, (.35,.35)*+[o][F-]{5}, (-.3,-.45)*+[o][F-]{6}
\end{xy}
\begin{xy}0;<.5cm,0cm>:
(-1,-1); (-1,1) **@{-}, (1,-1) **@{-}, (1,1); (1,-1) **@{-}, (-1,1) **@{-}, (-.15,-.1)*+[o][F-]{2}
\end{xy}
\begin{xy}0;<.5cm,0cm>:
(-1,-1); (-1,1) **@{-}, (1,-1) **@{-}, (1,1); (1,-1) **@{-}, (-1,1) **@{-}
\end{xy}
\begin{xy}0;<.5cm,0cm>:
(-1,-1); (-1,1) **@{-}, (1,-1) **@{-}, (1,1); (1,-1) **@{-}, (-1,1) **@{-}, (-.45,.25)*+[o][F-]{1}, (.4,-.2)*+[o][F-]{4}
\end{xy}
\begin{xy}0;<.5cm,0cm>:
(-1,-1); (-1,1) **@{-}, (1,-1) **@{-}, (1,1); (1,-1) **@{-}, (-1,1) **@{-}
\end{xy}
\begin{xy}0;<.5cm,0cm>:
(-1,-1); (-1,1) **@{-}, (1,-1) **@{-}, (1,1); (1,-1) **@{-}, (-1,1) **@{-}
\end{xy}
\begin{xy}0;<.5cm,0cm>:
(-1,-1); (-1,1) **@{-}, (1,-1) **@{-}, (1,1); (1,-1) **@{-}, (-1,1) **@{-}, (0,0)*+[o][F-]{3}
\end{xy}\]
We will use a shorthand notation for pictures like the one above, e.g., \[ ||56|2||14|||3 \] Here we draw a vertical bar for the divisions between the boxes, so that the number of boxes is one more than the number of bars. Also notice that by way of standardization, we list the balls in each box in increasing order. So balls 5 and 6 were placed in the third box, ball 2 in the fourth box, and so on. We call such diagram for balls in boxes a \emph{barred permutation}.

We can partition the set of all barred permutations (arrangements of balls in boxes) according to the underlying permutations. Since we require numbers to increase within a box, we know that there \emph{must} be a bar in each descent position of a barred permutation (i.e., at the end of every maximal increasing run), but otherwise, we can insert bars into the gaps between the numbers at will. So, for example, the barred permutations corresponding to $w = 562143$ are obtained from the barred permutation
\[ 56|2|14|3 \] by adding bars arbitrarily in the gaps between the numbers, i.e., there can be any number of bars to the left of 5, any number of bars between the 5 and the 6, at least one bar between the 6 and the 2, at least one bar between the 2 and the 1, and so on.

With the identity permutation, $w = 12\cdots n$, there are no required bars and we are free to add any number of bars in each of the $n+1$ gaps between the numbers. Thus, the generating function for barred permutations corresponding to  $12\cdots n$  is the generating function for the ``multi-choose" numbers $\left(\binom{n+1}{k}\right) = \binom{k+n}{k} = \binom{k+n}{n}$, i.e., the number of ways to choose, with repetition, $k$ things from an $(n+1)$-element set. With the convention that each barred permutation has weight $t^{\# \mbox{\tiny of bars }+1}$ ($=t^{\# \mbox{\tiny of boxes}}$), we express this generating function as:
\begin{align*}
\sum_{k\geq 0}\binom{k+n-1}{n}t^k &=t\cdot\underbrace{(1+t+t^2+\cdots)}_{\mbox{\tiny bars in first gap}}\underbrace{(1+t+t^2+\cdots)}_{\mbox{\tiny bars in second gap}}\cdots \\
 &= t\cdot \underbrace{(1+t+t^2+\cdots)^{n+1}}_{\mbox{\tiny $n+1$ gaps}} \\
 &= \frac{t}{(1-t)^{n+1}}.
\end{align*}

For permutations with descents, we need only multiply this generating function by a power of $t$ to reflect the number of bars required by descent positions: 
\begin{align}
\frac{t^{\des(w)+1}}{(1-t)^{n+1}} &= t^{\des(w)}\cdot \sum_{k\geq 0} \binom{k+n-1}{n} t^k \nonumber\\
&=\sum_{k\geq 0} \binom{k+n-1}{n} t^{\des(w)+k}\nonumber \\
&= \sum_{l\geq 0} \binom{l+n-1-\des(w)}{n} t^l.\label{eq:wsum}
\end{align}
So, returning to the example of $w=562143$, we get that the generating function for the sequence counting its corresponding barred permutations according to the number of bars is: \[\frac{t^4}{(1-t)^7}= t^3\cdot \sum_{k\geq 0} \binom{k+5}{6}t^k = \sum_{l\geq 0} \binom{l+2}{6} t^l. \]

Summing \eqref{eq:wsum} over all permutations in $S_n$ gives us, on the one hand, \[ \frac{A_n(t)}{(1-t)^{n+1}}.\] On the other hand, the sum gives the generating function for all barred permutations, which we know to be $\sum_{k\geq 0} k^n t^k$.

In other words, we have the following.

\begin{thm}[The generating function]\label{thm:Eul}
For any $n\geq 1$, we have
\[ \frac{A_n(t)}{(1-t)^{n+1}} = \sum_{k \geq 0} k^n t^k.\]
\end{thm}

This theorem is a classic result in enumerative combinatorics, going back at least to work of MacMahon \cite[Chapter IV, \S 462]{M}. See \cite[Section 4.5]{St} for a modern treatment and generalizations.

So, for example, the reader can check the $n=3$ case gives the generating function for the sequence $1, 8, 27, \ldots$, of cubes: \[ \frac{t + 4t^2 + t^3}{(1-t)^4} = t + 8t^2 + 27t^3 + 64 t^4 + 125t^5 +\cdots.\] 

Another consequence of the balls in boxes approach is \emph{Worpitzky's identity}, which is  explained in Knuth's book \cite[Section 5.1.3]{K}:
\begin{equation}\label{eq:worp}
 k^n = \sum_{1\leq i \leq n} A_{n,i} \binom{k+n-i}{n}.
\end{equation}
This equation shows Eulerian numbers describe a kind of ``change-of-basis" between binomial coefficients and $n$th powers.

For example, when $n=4$, $k=3$, \[ A_{4,1}\binom{6}{4} + A_{4,2}\binom{5}{4} + A_{4,3}\binom{4}{4} + A_{4,4}\binom{3}{4} = 1\cdot 15 + 11\cdot 5 + 11\cdot 1 + 1\cdot 0 = 81.\] 

To obtain Worpitzky's identity, we first recall from Equation \eqref{eq:wsum} that \[ \frac{t^i}{(1-t)^{n+1}} = \sum_{k\geq 0} \binom{k+n-i}{n} t^k.\] Then from the left-hand side of Theorem \ref{thm:Eul}, we get: \begin{align*} \frac{A_n(t)}{(1-t)^{n+1}} &= \sum_{i=0}^{n} A_{n,i}\left(\frac{t^i}{(1-t)^{n+1}}\right)\\
&= \sum_{i=0}^n A_{n,i}\cdot\sum_{k\geq 0} \binom{k+n-i}{n}t^k\\
&=\sum_{k\geq 0} \left(\sum_{i=0}^n A_{n,i}\binom{k+n-i}{n}\right) t^k.
\end{align*}
But according to the right-hand side of Theorem \ref{thm:Eul}, the coefficient of $t^k$ for this series is $k^n$. Hence, \eqref{eq:worp} follows.

Next, we can derive a recurrence for Eulerian polynomials with Theorem \ref{thm:Eul} and some calculus. Let $F_n(t) = \sum k^n t^k = A_n(t)/(1-t)^{n+1}$, and observe that 
\begin{align*}
F_n(t) &= \sum_{k\geq0} k^n t^k \\
  &= t\cdot \sum_{k\geq 0} k^{n-1} \cdot k t^{k-1} \\
  &= t\cdot F_{n-1}'(t) \\
  &= t\left( \frac{nA_{n-1}(t)}{(1-t)^{n+1}} + \frac{A'_{n-1}(t)}{(1-t)^n} \right)\\
  &= \frac{ntA_{n-1}(t) + t(1-t)A_{n-1}'(t)}{(1-t)^{n+1}}.
\end{align*}
Comparing numerators yields the following well known identity for Eulerian polynomials: 
\begin{equation}\label{eq:Apoly}
A_n(t) = ntA_{n-1}(t) + t(1-t)A'_{n-1}(t).
\end{equation}
From this equation, it is possible to deduce that the Eulerian polynomials factor completely over the reals, and moreover, the roots of $A_n(t)$ are all distinct, nonpositive numbers. According to Carlitz \cite{Car} this fact was first observed by Frobenius.

By comparing the coefficient of $t^i$ on the left and on the right of \eqref{eq:Apoly}, we also get a handy  recurrence for Eulerian numbers:
\begin{equation}\label{eq:Arec}
A_{n,i} = iA_{n-1,i} + (n+1-i)A_{n-1,i-1},
\end{equation}
which resembles a weighted version of Pascal's recurrence for binomial coefficients. 

The reader should check this numeric recurrence  against Table \ref{tab:Eul} to see how nicely it works. When drawing the numbers in a triangular array, we don't simply add the two numbers above to get the next entry, we add the appropriately weighted linear combination of the two. For example, the first few rows of such a triangle can be generated by hand as follows (so no need to write out all 120 permutations in $S_5$ to get the numbers $A_{5,i}$):
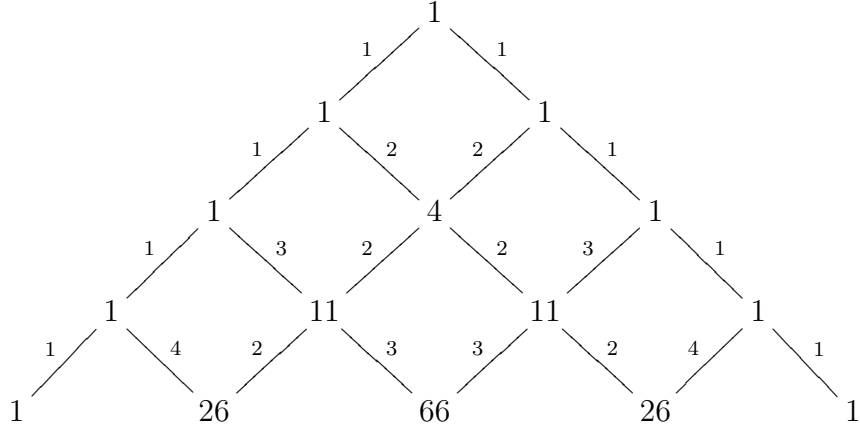
\[ \xymatrix{
    & & & & 1 \ar@{-}[dl]_{1}
       \ar@{-}[dr]^{1}& \\
    & & & 1 \ar@{-}[dl]_{1}
       \ar@{-}[dr]^{2}& & 1 \ar@{-}[dl]_{2}
       \ar@{-}[dr]^{1}& \\
    & & 1 \ar@{-}[dl]_{1}
       \ar@{-}[dr]^{3}& & 4 \ar@{-}[dl]_{2}
       \ar@{-}[dr]^{2} & & 1 \ar@{-}[dl]_{3}
       \ar@{-}[dr]^{1} \\
    & 1 \ar@{-}[dl]_{1}
       \ar@{-}[dr]^{4} & & 11 \ar@{-}[dl]_{2} 
       \ar@{-}[dr]^{3} & & 11 \ar@{-}[dl]_{3} 
       \ar@{-}[dr]^{2}& & 1 \ar@{-}[dl]_{4} 
       \ar@{-}[dr]^{1}\\
     1 & & 26 & & 66 & & 26 & & 1  }
\]

There is also a direct combinatorial proof of this recurrence. Imagine inserting the number $n$ into a permutation of $S_{n-1}$. What happens if you insert it in a descent position? What happens if you insert it in an ascent position? See, for example, \cite[Section 5.1.3]{K}.

There are many, many, more things that could be said at this point about Eulerian numbers. Some of these, including the row symmetry an astute reader may have noticed in Table \ref{tab:Eul}, will be discussed at the end of the paper.

\subsection*{Two-sided Eulerian numbers}

We are now going to extend the method of balls in boxes to study the ``two-sided" Eulerian polynomial, \[ A_n(s,t) := \sum_{w \in S_n} s^{\des(w^{-1})+1}t^{\des(w)+1} = \sum_{i,j=1}^n A_{n,i,j} s^i t^j.\] We think of $A_n(s,t)$ as the generating function for the joint distribution of left and right descents since $A_{n,i,j}$ is the number of permutations $w$ with $i-1$ left descents and $j-1$ right descents, i.e., \[A_{n,i,j} = |\{ w \in S_n : \des(w^{-1})=i-1 \mbox{ and } \des(w) = j-1\}|.\] In Table \ref{tab:Eul2} we see the two-sided Eulerian numbers for $n\leq 8$.

The two-sided Eulerian numbers were first studied by Carlitz, Roselle, and Scoville in 1966 \cite{CRS}, though rather than descents and inverse descents, they looked at the equivalent notion of ``jumps" (ascents) and ``readings" (inverse ascents). The results we prove here are all presented in \cite[Section 7]{CRS}, and proved using a mixture of combinatorics and manipulatorics (i.e., manipulation of formulas using binomial identities and such). Here we take the balls in boxes approach to prove:
\begin{itemize}
\item an analogue of Theorem \ref{thm:Eul}, 
\item a Worpitzky-like identity to compare with \eqref{eq:worp}, and 
\item a recurrence relation that refines \eqref{eq:Arec}.
\end{itemize}

\begin{table}
\begin{center}
{\footnotesize
\[ \begin{array}{c}  
  n=1\\
  \begin{array}{r|c}
  i\backslash j&1\\
  \hline
  1&1 \end{array}
  \end{array}  
\quad   
\begin{array}{c}
n=2 \\
\begin{array}{r| c c}
 i\backslash j & 1 & 2 \\
  \hline
 1 & 1 & 0 \\
 2 &0 & 1 
\end{array}
\end{array}
\quad
\begin{array}{c}
n=3 \\
\begin{array}{r|c c c}
i\backslash j & 1 & 2 & 3\\
\hline
1 & 1 & 0 & 0 \\
2& 0 & 4 & 0 \\
3 & 0 & 0 & 1 \\ 
\end{array}
\end{array}
\quad
\begin{array}{c}
n=4\\
\begin{array}{r|c c c c}
i\backslash j & 1 & 2 & 3 & 4\\
\hline 
1& 1 & 0 & 0 & 0\\
2& 0 & 10 & 1 & 0 \\
3& 0 & 1 & 10 & 0 \\ 
4& 0 & 0 & 0 & 1
\end{array} \\
\end{array}
\]
\[
\begin{array}{c}
n=5\\
\begin{array}{r|c c c c c}
i \backslash j & 1 & 2 & 3 & 4 & 5\\
\hline
1& 1 & 0 & 0 & 0 & 0\\
2& 0 & 20 & 6 & 0 & 0\\
3& 0 & 6 & 54 & 6 & 0 \\ 
4& 0 & 0 & 6 & 20 & 0\\
5& 0 & 0 & 0 & 0 & 1
\end{array}
\end{array}
\quad
\begin{array}{c}
n=6\\
\begin{array}{r|c c c c c c}
i \backslash j & 1 & 2 & 3 & 4 & 5 & 6\\
\hline
1& 1 & 0 & 0 & 0 & 0 &0\\
2& 0 & 35 & 21 & 1 & 0 &0\\
3& 0 & 21 & 210 & 70 & 1 & 0 \\ 
4& 0 & 1 & 70 & 210 & 21 & 0\\
5& 0 & 0 & 1 & 21 & 35 & 0\\
6& 0 & 0 & 0 & 0 & 0 & 1
\end{array}\end{array}
\]
\[
\begin{array}{c}
n=7\\
\begin{array}{r|c c c c c c c}
i \backslash j & 1 & 2 & 3 & 4 & 5 & 6 & 7\\
\hline
1& 1 & 0 & 0 & 0 & 0 &0 & 0\\
2& 0 & 56 & 56 & 8 & 0 &0 & 0\\
3& 0 & 56 & 659 & 440 & 36 & 0 & 0\\ 
4& 0 & 8 & 440 & 1520 & 440 & 8 & 0\\
5& 0 & 0 & 36 & 440 & 659 & 56 & 0\\
6& 0 & 0 & 0 & 8 & 56 & 56 & 0\\
7& 0 & 0 & 0 & 0 & 0 & 0& 1
\end{array}
\end{array}
\]
\[
\begin{array}{c}
n=8\\
\begin{array}{r|c c c c c c c c}
i\backslash j & 1 & 2 & 3 & 4 & 5 & 6 & 7 & 8\\
\hline
1 & 1 & 0 & 0 & 0 & 0 &0 & 0 & 0\\
2& 0 & 84 & 126 & 36 & 1 &0 & 0& 0\\
3& 0 & 126 & 1773 & 1980 & 405 & 9 &0 & 0\\ 
4& 0 & 36 & 1980 & 8436 & 4761 & 405 & 1& 0\\
5& 0 & 1 & 405 & 4761 & 8436 & 1980 & 36& 0\\
6& 0 & 0 & 9 & 405 & 1980 & 1773 & 126 & 0\\
7& 0 & 0 & 0 & 1 & 36 & 126 & 84 & 0 \\
8& 0 & 0 & 0 & 0 & 0 &0 & 0 & 1
\end{array}
\end{array}
\]
}
\medskip
\caption{The two-sided Eulerian numbers $[A_{n,i,j}]_{1\leq i,j\leq n}$, $n=1,\ldots,8$ }\label{tab:Eul2}
\end{center}
\end{table}

\subsection*{Balls in boxes, revisited}

To prove these ``two-sided" results, we will put balls in boxes another way. First, we introduce a two-dimensional analogue of barred permutations. 

Permutations in $S_n$ can be represented visually as an array of $n$ indistinguishable balls so that no two balls lie in the same column or row. If $w(i) = j$, we put a ball in column $i$ (from left to right) and row $j$ (from bottom to top). For example, with $w = 562143$, we draw:
\[ 
\begin{array}{r| c c c c c c }
6 & & \bullet & & & & \\
5 & \bullet & & & & & \\
4 & & & & & \bullet & \\
3 & & & & & & \bullet \\
2 & & & \bullet & & & \\
1 & & & & \bullet & & \\
\hline
j/i & 1 & 2 & 3 & 4 & 5 & 6 \\
\end{array}
\]

The advantage of such an array is that $w^{-1}$  is easily seen in this picture. Rather than reading the heights of the balls from left to right to get $w$, we can read the column numbers of the balls from bottom to top to get $w^{-1}$. So in the example above, $w^{-1} = 436512$. (Check this! Why is it generally true?)

The analogue of a barred permutation, which we call a \emph{two-sided} barred permutation, is any way of inserting both horizontal and vertical lines into the array of balls, with the requirement that there be at least one vertical line between balls that form a descent in $w$, and at least one horizontal line between balls that form a descent in $w^{-1}$. Other horizontal and vertical lines can be added arbitrarily. For example:
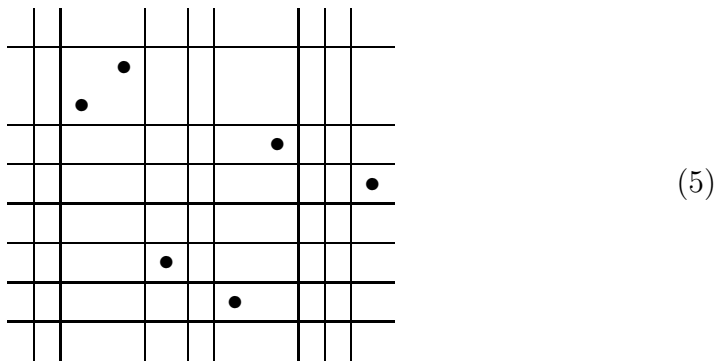
\begin{equation}\label{eq:2bar}
\begin{array}{c | c | c c | c  | c | c c | c| c | c }
 & & & & & & & & & \\
\hline
 & & & \bullet & & & & & & &\\
 & & \bullet & & & & & & & & \\
\hline
 & & & & & & & \bullet & & & \\
\hline
 & & & & & & & & & & \bullet \\
\hline
 & & & & & & & & &\\
\hline
 & & & & \bullet & & & & & & \\
\hline
 & & & & & & \bullet & & & & \\
\hline
 & & & & & & & & &
\end{array}
\end{equation}
In a sense, this two-sided barred permutation corresponds to two ordinary barred permutations, one for $w$: $||56|2||14|||3$, and one for $w^{-1}$: $|4|3||6|5|12|$.

Now fix an arrangement of balls corresponding to a permutation $w$. We will give a two-sided barred permutation the weight \[ s^{(\# \mbox{\tiny of horizontal bars } +1)}t^{(\# \mbox{\tiny of vertical bars } +1)},\] so the example seen in \eqref{eq:2bar} would contribute $s^8t^9$. Since the horizontal and vertical bars can be inserted independently of one another, we see that the generating function for the number of two-sided barred permutations corresponding to a fixed permutation $w$ in $S_n$ is the product of the generating function for barred permutations for $w^{-1}$ (in the variable $s$) with the generating function for barred permutations of $w$ (in the variable $t$). Thus \eqref{eq:wsum} gives: \begin{equation}\label{eq:w}
 \frac{s^{\des(w^{-1})+1} t^{\des(w) +1}}{(1-s)^{n+1}(1-t)^{n+1}}.
\end{equation}
Adding up \eqref{eq:w} over all permutations in $S_n$ we get the generating function for all two-sided barred permutations is \begin{equation}\label{eq:twosum}
 \frac{ A_n(s,t)}{(1-s)^{n+1}(1-t)^{n+1}}.
\end{equation}

Now consider forming two-sided barred permutations with the bars first. Given $k-1$ vertical bars and $l-1$ horizontal bars, we get a $k$-by-$l$ grid of boxes in which to place our balls, with the convention that if more than one ball goes into a particular row or column, we arrange the balls diagonally from bottom left to top right. For example, the following arrangement of 7 balls in a $5\times 4$ grid of boxes yields the following two-sided barred permutation (with underlying permutation $1723465$): 
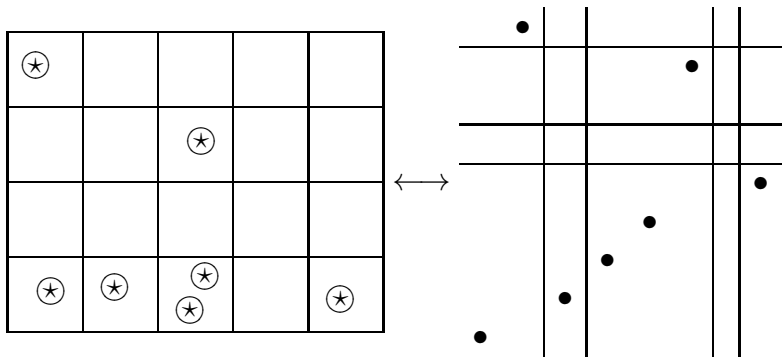
\[
\begin{xy}0;<.5cm,0cm>:
(-1,-4); (-1,4) **@{-}, (1,-4) **@{-}, (1,4); (1,-4) **@{-}, (-1,4) **@{-}, (-1,-2); (1,-2) **@{-}, (-1,0); (1,0) **@{-}, (-1,2); (1,2) **@{-}, (-.25,3.1)*+[o][F-]{\star}, (.15,-2.9)*+[o][F-]{\star}
\end{xy}
\begin{xy}0;<.5cm,0cm>:
(-1,-4); (-1,4) **@{-}, (1,-4) **@{-}, (1,4); (1,-4) **@{-}, (-1,4) **@{-}, (-1,-2); (1,-2) **@{-}, (-1,0); (1,0) **@{-}, (-1,2); (1,2) **@{-}, (-.15,-2.8)*+[o][F-]{\star}
\end{xy}
\begin{xy}0;<.5cm,0cm>:
(-1,-4); (-1,4) **@{-}, (1,-4) **@{-}, (1,4); (1,-4) **@{-}, (-1,4) **@{-}, (-1,-2); (1,-2) **@{-}, (-1,0); (1,0) **@{-}, (-1,2); (1,2) **@{-}, (-.15,-3.4)*+[o][F-]{\star}, (.25,-2.5)*+[o][F-]{\star}, (.15,1.1)*+[o][F-]{\star}
\end{xy}
\begin{xy}0;<.5cm,0cm>:
(-1,-4); (-1,4) **@{-}, (1,-4) **@{-}, (1,4); (1,-4) **@{-}, (-1,4) **@{-}, (-1,-2); (1,-2) **@{-}, (-1,0); (1,0) **@{-}, (-1,2); (1,2) **@{-}, 
\end{xy}
\begin{xy}0;<.5cm,0cm>:
(-1,-4); (-1,4) **@{-}, (1,-4) **@{-}, (1,4); (1,-4) **@{-}, (-1,4) **@{-}, (-1,-2); (1,-2) **@{-}, (-1,0); (1,0) **@{-}, (-1,2); (1,2) **@{-},  (-.15,-3.1)*+[o][F-]{\star}
\end{xy}
\longleftrightarrow
\begin{array}{c c | c | c c c | c  | c   }
 & \bullet & & & & & &  \\
\hline
 & & & & & \bullet & & \\
 & &  & & & & & \\
\hline
 & & & & & & & \\
\hline
 & & & & & & &  \bullet \\
 & & & &\bullet & & &\\
 & & & \bullet  & & & & \\
 & & \bullet & & & & & \\
 \bullet & & & & & & &
\end{array}
\]

With $n$ unlabeled balls and $kl$ distinct boxes, this means there are a total of \[ \left(\binom{kl}{n}\right) = \binom{ kl +n -1 }{n} \] two-sided barred permutations with $k-1$ vertical lines and $l-1$ horizontal lines. (We are essentially choosing $n$ of the $kl$ boxes, with repetition allowed.) The generating function for the number of \emph{all} two-sided barred permuations is thus \[ \sum_{k,l \geq 0} \binom{kl+n-1}{n} s^k t^l,\]  and comparing with \eqref{eq:twosum} yields the following.

\begin{thm}[The two-sided generating function]\label{thm:Eul2}
For $n\geq 1$, we have
\[ \frac{A_n(s,t)}{(1-s)^{n+1}(1-t)^{n+1}} = \sum_{k,l\geq0} \binom{kl+n-1}{n} s^k t^l.\]
\end{thm}

Just as Worpitzky's identity can be derived from Theorem \ref{thm:Eul} and Equation \eqref{eq:wsum}, Theorem \ref{thm:Eul2} and Equation \eqref{eq:wsum} yields the following ``Worpitzky-like" identity of binomial coefficients.
\begin{equation}\label{eq:worp2} \binom{kl+n-1}{n} = \sum_{i,j=1}^n A_{n,i,j} \binom{k+n-i}{n} \binom{l+n-j}{n}.
\end{equation}

To get \eqref{eq:worp2}, we use \eqref{eq:wsum} and \eqref{eq:w} to see that for any permutation $w \in S_n$ with $\des(w^{-1}) = i-1$, $\des(w) = j-1$, the generating function for its barred permutations is
\begin{align*}
 \frac{s^i t^j}{(1-s)^{n+1}(1-t)^{n+1}} &= \sum_{k\geq0} \binom{k+n-i}{n} s^k \cdot \sum_{l\geq0} \binom{l+n-j}{n} t^l \\
 &= \sum_{k,l\geq0} \binom{k+n-i}{n} \binom{l+n-j}{n} s^k t^l.
\end{align*}
There are $A_{n,i,j}$ such permutations, so summing over $w \in S_n$ yields
\begin{align*} \frac{A_n(s,t)}{(1-s)^{n+1}(1-t)^{n+1}} &= \sum_{i,j=1}^n A_{n,i,j} \left(\frac{s^i t^j}{(1-s)^{n+1}(1-t)^{n+1}} \right)\\
&= \sum_{i,j=1}^n A_{n,i,j}\cdot\sum_{k,l\geq 0}\binom{k+n-i}{n}\binom{l+n-j}{n} s^kt^l \\
&= \sum_{k,l\geq 0} \left(\sum_{i,j=1}^n A_{n,i,j}\binom{k+n-i}{n}\binom{l+n-j}{n}\right) s^kt^l
\end{align*} 
as the generating function for all two-sided barred permutations for $S_n$. But we already know there are $\binom{kl+n-1}{n}$ two-sided barred permutations of weight $s^kt^l$, and so we have \eqref{eq:worp2}, as desired.

We remark that this Worpitzky-like identity can also be proved (and considerably generalized) using the method of bipartite $P$-partitions, as follows from \cite[Corollary 10]{Ges}.

Now from Theorem \ref{thm:Eul2} we obtain a recurrence relation akin to \eqref{eq:Apoly}, which makes for easy computation of the two-sided Eulerian polynomials. See \cite[Equation 7.8]{CRS}. We have:
\begin{align}
nA_n(s,t) &= \big( n^2st + (n-1)(1-s)(1-t) \big)A_{n-1}(s,t) \nonumber \\
&\quad + nst(1-s)\frac{\partial}{\partial s} A_{n-1}(s,t) + nst(1-t)\frac{\partial}{\partial t} A_{n-1}(s,t) \nonumber \\
&\quad + st(1-s)(1-t)\frac{\partial^2}{\partial s \partial t} A_{n-1}(s,t). \label{eq:Apoly2}
\end{align} 

To obtain \eqref{eq:Apoly2}, we let $F_n(s,t) = A_n(s,t)/(1-s)^{n+1}(1-t)^{n+1}$, and note the following identity of binomial coefficients that we will need: 
\begin{align*} 
n\binom{kl+n-1}{n} &=\frac{n\cdot (kl+n-1)!}{n!(kl-1)!} \\
&=(kl+n-1)\cdot \frac{(kl+n-2)!}{(n-1)!(kl-1)!} \\
&= kl\binom{kl+n-2}{n-1} + (n-1)\binom{kl+n-2}{n-1}.
\end{align*}

Thus, we can see that
\begin{align}
nF_n(s,t) &= \sum_{k,l\geq 0} n\binom{kl+n-1}{n} s^k t^l \nonumber \\
 &= \sum_{k,l \geq 0} kl\binom{kl+n-2}{n-1} s^k t^l + \sum_{k,l \geq 0} (n-1)\binom{kl+n-2}{n-1}s^k t^l  \nonumber \\
 &=st \frac{\partial^2}{\partial s \partial t} F_{n-1}(s,t) + (n-1)F_{n-1}(s,t).\label{eq:F}
\end{align}

Now, with a little calculus, we have
\begin{align*}
 \frac{\partial^2}{\partial s \partial t} F_{n-1}(s,t) &=  \frac{\partial}{\partial s}\left[ \frac{n A_{n-1}(s,t)}{(1-s)^n(1-t)^{n+1}} + \frac{\frac{\partial}{\partial t} A_{n-1}(s,t)}{(1-s)^n(1-t)^n} \right] \\ 
&= \frac{n^2 A_{n-1}(s,t)}{(1-s)^{n+1}(1-t)^{n+1}} + \frac{n \frac{\partial}{\partial s}A_{n-1}(s,t)}{(1-s)^n (1-t)^{n+1}} \\
& \quad  + \frac{n \frac{\partial}{\partial t}A_{n-1}(s,t)}{(1-s)^{n+1} (1-t)^{n}} + \frac{\frac{\partial^2}{\partial s \partial t} A_{n-1}(s,t)}{(1-s)^n(1-t)^n}.
\end{align*}

After giving each term a common denominator, we  compare numerators in \eqref{eq:F} to get
\eqref{eq:Apoly2},
as desired.

By comparing coefficients on both sides of \eqref{eq:Apoly2}, we also get the following four-term numeric recurrence:
\begin{align}
nA_{n,i,j} &=(ij+n-1)A_{n-1,i,j} \nonumber \\
& \quad + [1-n+j(n+1-i)]A_{n-1,i-1,j} + [1-n+i(n+1-j)]A_{n-1,i,j-1} \nonumber \\
& \quad + [n-1+(n+1-i)(n+1-j)]A_{n-1,i-1,j-1}. \label{eq:Arec2}
\end{align}

The reader is invited to use this recurrence to obtain the first few arrays in Table \ref{tab:Eul2}. Is there a nice visual way to understand this recurrence?

\subsection*{Valley-hopping and $\gamma$-nonnegativity}

A sharp-eyed reader may have already noticed the well-known fact that the Eulerian numbers are \emph{symmetric} for fixed $n$: \[A_{n,i} = A_{n,n+1-i},\] or in terms of Eulerian polynomials, $A_n(t) = t^{n+1}A_n(1/t)$.  This symmetry can be explained by observing that if $w$ has $i$ descents it has $n-1-i$ ascents.

Similarly, one can find the following symmetries for two-sided Eulerian numbers in Table \ref{tab:Eul2}:
\begin{align}
A_{n,i,j} = A_{n,j,i}, \quad &\mbox{or} \quad A_n(s,t) = A_n(t,s), \label{eq:sym1} \\
A_{n,i,j} = A_{n,n+1-i,n+1-j}, \quad &\mbox{or}\quad A_n(s,t) = (st)^{n+1}A_n(1/s,1/t), \mbox{ and} \label{eq:sym2} \\
A_{n,i,j} = A_{n,n+1-j,n+1-i}, \quad &\mbox{or} \quad A_n(s,t) = (st)^{n+1}A_n(1/t,1/s). \label{eq:sym3}
\end{align}

It is a fun exercise to come up with combinatorial arguments to verify symmetries \eqref{eq:sym1} and \eqref{eq:sym2}. (Hint: how does flipping a permutation array upside-down affect descents and inverse descents?) Symmetry \eqref{eq:sym3} follows from the first two.

Another property exhibited by the Eulerian polynomials is \emph{unimodality}, i.e., the Eulerian numbers in a given row increase up to a certain maximum and then decrease: \[ A_{n,1} \leq A_{n,2} \leq \cdots \leq A_{n,\lceil n/2 \rceil} \geq \cdots \geq A_{n,n-1} \geq A_{n,n}.\] Most good distributions satisfy this property of having the bulk of the mass is in the middle and the rest is spread out symmetrically. The canonical example is the binomial distribution, $\binom{n}{i}$ with fixed $n$. Similarly, the two-sided Eulerian numbers appear to increase in the direction of the middle of the main diagonal, with a maximum at $A_{n, \lceil n/2 \rceil, \lceil n/2 \rceil}$. However, things are somewhat delicate. For example, when $i< j \leq n/2$ it is not always true that $A_{n,i+1,j} \geq A_{n,i,j} \leq A_{n,i,j-1}$, i.e., moving toward the main diagonal can sometimes lead to a smaller number. The first examples of this occur for $n=8$, e.g., $A_{8,2,3} = 126 > 84= A_{8,2,2}$, and $A_{8,3,4}=1980>1773 = A_{8,3,3}$. See Table \ref{tab:Eul2}.

Both the symmetry and unimodality of the Eulerian polynomials follow from a stronger result, first proved in 1970 by Foata and Sch\"utzenberger \cite[Th\'eor\`eme 5.6]{FS}. See also Carlitz and Scoville \cite{Car2}.

\begin{thm}\label{thm:gamma}
For $n\geq 1$, there exist nonnegative integers $\gamma_{n,i}$, $i=1,\ldots,\lceil n/2\rceil$, such that
\[ A_n(t) = \sum_{1\leq i \leq \lceil n/2\rceil} \gamma_{n,i} t^i (1+t)^{n+1-2i}.\]
\end{thm}

For example, when $n=4,5$ we have: 
\begin{align*}
 A_4(t) &= t + 11t^2+11t^3 + t^4 \\
  &=  t(1+t)^3 + 8t^2(1+t),\\
A_5(t) &= t+ 26t^2 + 66t^3 + 26t^4 + t^5\\
 &= t(1+t)^4 + 22t^2(1+t)^2 + 16t^3.
 \end{align*}
In other words, the polynomials $A_n(t)$ can be expressed as a positive sum of symmetric binomial terms with the same center of symmetry. Since the binomial distribution is symmetric and unimodal, so is the Eulerian distribution. One might say the Eulerian distribution is ``super binomial". Gessel has conjectured a similar expansion for two-sided Eulerian numbers that we will discuss shortly. 

First, we explain how Foata and Sch\"utzenberger's result can be given a wonderfully visual proof via an action called \emph{valley-hopping}. The argument here has its roots in work of Foata and Strehl \cite{FSt} from 1974, was re-discovered by Shapiro, Woan, and Getu \cite{SWG} in 1983, and was dusted off more recently (and applied in a wonderful variety of ways) by Br\"and\'en \cite{Br} in 2008. 

We define the valley-hopping action on permutations as follows. We visualize permutations as arrays of balls again, but now we connect the dots to form a kind of mountain range. Some balls sit at peaks, others sit in valleys, and the rest are somewhere in between. If a ball is not at a peak or in a valley, it is free to jump straight across a valley to the nearest point on a slope at the same height.

Valley hopping naturally partitions $S_n$ into equivalence classes according to whether one permutation can be obtained from another through a sequence of hops. For example, the permutation $w=863247159$ would be drawn as follows:
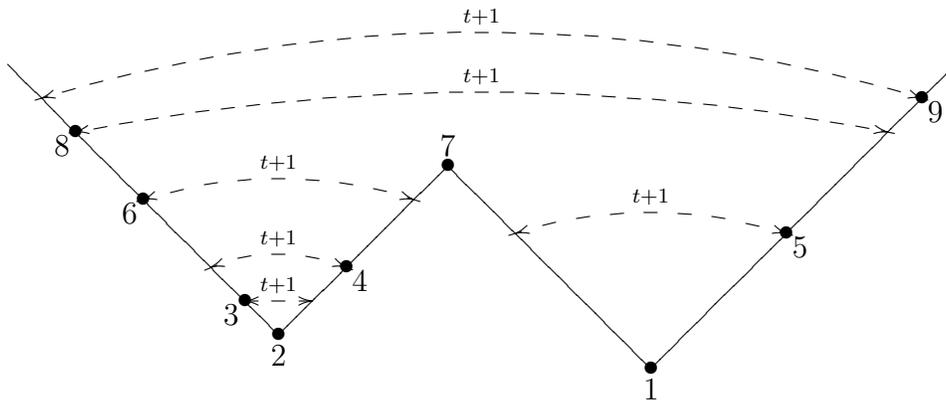
\begin{equation}\label{eq:val}
\begin{xy}0;<.9cm,0cm>:
(-2.5,-1.5); (-6.5,2.5) **@{-}, (0,1) **@{-}, (3,-2); (0,1) **@{-}, (7.5,2.5) **@{-},
(-2.5,-1.5)*{\bullet}, (-2.5,-1.8)*{2}, (0,1)*{\bullet}, (0,1.3)*{7}, (3,-2)*{\bullet}, (3,-2.3)*{1}, (7,2)*{\bullet}, (7.2,1.8)*{9}, (-3,-1)*{\bullet}, (-3.2,-1.2)*{3}, (-1.5,-.5)*{\bullet}, (-1.3,-.7)*{4}, (5,0)*{\bullet}, (5.2,-.2)*{5}, (-4.5,.5)*{\bullet}, (-4.7,.3)*{6}, (-5.5,1.5)*{\bullet}, (-5.7,1.3)*{8}, {\ar @{<-->}^{t+1} @/^15pt/ (-5.5,1.5); (6.5,1.5) }, {\ar @{<-->}^{t+1} @/^25pt/ (-6,2); (7,2) }, {\ar @{<-->}^{t+1} (-3,-1); (-2,-1) }, {\ar @{<-->}^{t+1} @/^5pt/ (-3.5,-.5); (-1.5,-.5) }, {\ar @{<-->}^{t+1} @/^8pt/ (-4.5,.5); (-.5,.5) }, {\ar @{<-->}^{t+1} @/^8pt/ (1,0); (5,0) }
\end{xy}
\end{equation}
There are $2^6$ permutations in its equivalence class, formed by choosing which of the six free balls will be on the left sides of their respective valleys and which will be on the right. Write $u\sim w$ if $u$ can be obtained from $w$ through valley-hopping. Notice that when a free ball is on the right side of a valley, it is \emph{not} in a descent position, while if it is on the left side of a valley, it \emph{is} in a descent position. (This is why we labeled the arcs in the picture with $t+1$.) Moreover, this property holds true regardless of the positions of the other free balls.

As for the non-free balls, we know that peaks are \emph{always} in descent positions while valleys are \emph{never} in descent positions. If a permutation has $i-1$ peaks, then it must have $i$ valleys, and the remaining $n+1-2i$ balls are free. Thus, we can conclude that for a fixed $w$ in $S_n$ with $i-1$ peaks, \[ \sum_{u \sim w} t^{\des(u)+1} = t^i(1+t)^{n+1-2i}.\]
For example, the equivalence class for $w=863247159$ would contribute \[ t^2(1+t)^6. \]

Since the union of all equivalence classes is $S_n$, we see that the Eulerian polynomial is a sum of terms of the form $t^i(1+t)^{n+1-2i}$, proving Theorem \ref{thm:gamma}. Moreover, the coefficient $\gamma_{n,i}$ equals the number of distinct equivalence classes with $i-1$ peaks.

We now turn to the two-sided case.

It is possible to show that any polynomial in two variables (of degree $n$ in each) satisfying symmetries \eqref{eq:sym1} and \eqref{eq:sym2} can be written uniquely in the basis \[ \{ (st)^i (s+t)^j (1+st)^{n+1-j-2i}\}_{0\leq j+2i \leq n+1}.\]
Thus, the two-sided Eulerian polynomials can be expressed in this basis, and Gessel has conjectured that such an expression is nonnegative. See \cite[Conjecture 10.2]{Br}. 

\begin{conj}[Gessel's conjecture]\label{conj:ges}
For $n\geq 1$, there exist nonnegative integers $\gamma_{n,i,j}$ , $0\leq i-1, j$, $j+2i \leq n+1$, such that \[ A_n(s,t) = \sum_{i,j} \gamma_{n,i,j} (st)^i (s+t)^j(1+st)^{n+1-j-2i}.\]
\end{conj}

For example, when $n=4,5$, we have
\begin{align*}
 A_4(s,t) &= st + 10(st)^2 + 10(st)^3 + (st)^4 + s^2t^3 + s^3t^2\\
 &= st(1+st)^3 + 7(st)^2(1+st) + (st)^2(s+t)\\
A_5(s,t) &= st + 20(st)^2 + 54(st)^3 + 20(st)^4 + (st)^5 + 6s^2t^3 + 6s^3t^2 + 6s^3t^4 + 6s^4t^3\\
&= st(1+st)^4 + 16(st)^2(1+st)^2 + 16(st)^3 + 6(st)^2(s+t)(1+st)
\end{align*}
In terms of the arrays $[A_{n,i,j}]$, we see:
\[ \left[\begin{array}{c c c c}
1 & 0 & 0 & 0\\
0 & 10 & 1 & 0 \\
0 & 1 & 10 & 0 \\ 
0 & 0 & 0 & 1
\end{array}\right] = \left[\begin{array}{c c c c}
1 & 0 & 0 & 0\\
0 & 3 & 0 & 0 \\
0 & 0 & 3 & 0 \\ 
0 & 0 & 0 & 1
\end{array}\right]
+
7\cdot\left[\begin{array}{c c c c}
0 & 0 & 0 & 0\\
0 & 1 & 0 & 0 \\
0 & 0 & 1 & 0 \\ 
0 & 0 & 0 & 0
\end{array}\right]
+\left[\begin{array}{c c c c}
0 & 0 & 0 & 0\\
0 & 0 & 1 & 0 \\
0 & 1 & 0 & 0 \\ 
0 & 0 & 0 & 0
\end{array}\right]
\]
and
\begin{align*}
\left[\begin{array}{c c c c c}
1 & 0 & 0 & 0 & 0\\
0 & 20 & 6 & 0 & 0\\
0 & 6 & 54 & 6 & 0 \\ 
0 & 0 & 6 & 20 & 0\\
0 & 0 & 0 & 0 & 1
\end{array}\right] &=
\left[\begin{array}{c c c c c}
1 & 0 & 0 & 0 & 0\\
0 & 4 & 0 & 0 & 0\\
0 & 0 & 6 & 0 & 0 \\ 
0 & 0 & 0 & 4 & 0\\
0 & 0 & 0 & 0 & 1
\end{array}\right]
+
16\cdot\left[\begin{array}{c c c c c}
0 & 0 & 0 & 0 & 0\\
0 & 1 & 0 & 0 & 0\\
0 & 0 & 2 & 0 & 0 \\ 
0 & 0 & 0 & 1 & 0\\
0 & 0 & 0 & 0 & 0
\end{array}\right]
\\
&+ 16\cdot\left[\begin{array}{c c c c c}
0 & 0 & 0 & 0 & 0\\
0 & 0 & 0 & 0 & 0\\
0 & 0 & 1 & 0 & 0 \\ 
0 & 0 & 0 & 0 & 0\\
0 & 0 & 0 & 0 & 0
\end{array}\right]
+ 6 \cdot\left[\begin{array}{c c c c c}
0 & 0 & 0 & 0 & 0\\
0 & 0 & 1 & 0 & 0\\
0 & 1 & 0 & 1 & 0 \\ 
0 & 0 & 1 & 0 & 0\\
0 & 0 & 0 & 0 & 0
\end{array}\right]
\end{align*}

It is possible that one might use a ``manipulatorics" approach to prove Conjecture \ref{conj:ges}. (Perhaps an inductive proof using the recurrence in \eqref{eq:Apoly2}?) 
However, a more satisfying proof might be one that generalizes the valley-hopping proof of Theorem \ref{thm:gamma}.

If we look at both descents and inverse descents for the valley-hopping equivalence class of $w=12\cdots n$ (i.e., the class with no peaks), we get a distribution of $st(1+st)^{n-1}$. This is encouraging, but for the class of $w=863247159$ shown in \eqref{eq:val}, we get \[ \sum_{u\sim w} s^{\des(u^{-1})+1}t^{\des(u)+1} = s^3t^2(1+t)^2(1+st)^4,\]
which is not even symmetric in $s$ and $t$. So valley-hopping as presently done does not immediately give us a way to prove Gessel's conjecture.

How should we partition $S_n$ so that we get groupings whose distribution of descents and inverse descents is given by $(st)^i (s+t)^j (1+st)^{n+1-j-2i}$?

\medskip
\noindent\textbf{Generalization}

We finish by remarking that $S_n$ is an example of a \emph{finite reflection group}, or \emph{Coxeter group}. The notion of a descent can be generalized to any Coxeter group, and there is a ``Coxeter-Eulerian" polynomial that enjoys many of the same properties  of the classical Eulerian polynomial, including an analogue of Theorem \ref{thm:gamma}. Moreover, this polynomial has topological meaning, as the ``$h$-polynomial" of something called the Coxeter complex. See \cite{BjBr, DPS, Stem}.

The two-sided Eulerian polynomial generalizes to Coxeter groups as well, and seems to enjoy many of the same properties of $A_n(s,t)$. In particular, the analogue of Conjecture \ref{conj:ges} appears to hold in any finite Coxeter group. It would be interesting to have a general approach to the problem.

\subsection*{Acknowledgements}

I would like to thank Ira Gessel for teaching me how to put balls in boxes and for sharing his conjecture.

\end{document}